\input amstex
\documentstyle{amsppt}
\magnification=\magstep1
\hsize=5.2in
\vsize=6.8in

\topmatter

\centerline  {\bf ORBIT INEQUIVALENT ACTIONS FOR GROUPS }

\centerline {\bf  CONTAINING  A COPY OF $\Bbb F_2$} 

\vskip .1in

\centerline {\rm ADRIAN IOANA}
\vskip 0.05in

\abstract
We prove that if a countable  group $\Gamma$ contains a copy of $\Bbb F_2$, then it admits uncountably many non orbit equivalent actions.
\endabstract

\endtopmatter

\document
\vskip 0.1in

\centerline{\S0. {\bf Introduction.}}
\vskip 0.1in
Throughout this paper we consider free, ergodic, measure preserving (m.p.) actions $\Gamma\curvearrowright (X,\mu)$ of countable, discrete  groups $\Gamma$ on  standard probability spaces $(X,\mu)$.
Measurable group theory is roughly the study of such group actions  from the viewpoint of the induced orbit equivalence relation.
A basic question in measurable group theory is to find groups $\Gamma$ which admit many non-orbit equivalent actions (see the survey [Sh2]).
 In this respect, recall that two free, ergodic, m.p. actions  $\Gamma\curvearrowright (X,\mu)$ and $\Lambda\curvearrowright (Y,\nu)$ are said to be {\bf orbit equivalent (OE)} if they induce isomorphic equivalence relations, i.e. if there exists a measure space isomorphism $\theta: (X,\mu)\rightarrow (Y,\nu)$ such that $\theta(\Gamma x)=\Lambda \theta(x),$ for almost every $x\in X$.

The striking lack of rigidity manifested by amenable groups (any two free, ergodic m.p. actions of any two infinite amenable groups $\Gamma$ and $\Lambda$ are orbit equivalent--a result proved by Dye in the case $\Gamma$ and $\Lambda$ are abelian ([D]) and by Ornstein--Weiss in general ([OW])) implies that the above question is well-posed only for non-amenable groups.  
For  a non-amenable group $\Gamma$, it is known that $\Gamma$ admits at least two non-{\bf OE} actions ([CW],[Sc],[H]).
Moreover, recently, the following classes of non-amenable groups have been shown to admit uncountably many non-{\bf OE} actions: property (T) groups ([H]), free groups ([GP], see also [T],[I2]), weakly rigid groups ([P2]), non-amenable products of infinite  groups ([P4], see also  [MSh],[I1]) and mapping class groups ([Ki]).  This classes of group added to the already known ones ([BG], [Z], [GG]). 

In this paper, we prove that the same is true for a new, large class of non-amenable groups.

\proclaim {Main result} Let $\Gamma$ be a countable discrete group which contains a copy of the free group $\Bbb F_2$. Then 
 $\Gamma$ has uncountably many non-OE actions. 
 \vskip 0.03in
Moreover, any group $\Lambda$ which is measure equivalent to $\Gamma$ has uncountably many non-OE actions.
\endproclaim
\vskip 0.05in
Note that this result covers most non-amenable groups. The question of whether
 every non-amenable $\Gamma$ contains a copy of $\Bbb F_2$, known as von Neumann's problem, was open for a long time, until it was settled in the negative by Ol'shanskii ([O]).
 
To outline the proof of the main result, recall that if we view $\Bbb F_2$ as a finite index subgroup of SL$_2(\Bbb Z)$, then the pair $(\Bbb F_2\ltimes\Bbb Z^2,\Bbb Z^2)$ has the relative property (T) of Kazhdan-Margulis ([K],[M]). This fact implies that the induced m.p. action $\Bbb F_2\curvearrowright^{\alpha} \Bbb T^2=\hat{\Bbb Z^2}$ is rigid, in the sense of Popa ([P1]). The main idea of the proof is then to consider the class $\Cal F$ of actions $\Gamma\curvearrowright X$ for which the restriction $\Bbb F_2\curvearrowright X$ admits $\alpha$ as a quotient. 
Note that the rigidity of $\alpha$  has been successfully used before by Gaboriau and Popa to show that non-abelian free groups admit uncountably many non-OE actions ([GP]) .

Using a separability argument (in the spirit of [C2], [P1] and [GP]) in connection with the rigidity of $\alpha$, we prove that, modulo countable sets, orbit equivalence of two $\Gamma$-actions from $\Cal F$ implies conjugacy of the restrictions to $\Bbb F_2$ (Theorem 1.3.). 
On the other hand, using the co-induced construction (see Section 2) we provide uncountably many actions in $\Cal F$ for which the restrictions to $\Bbb F_2$ are mutually non-conjugate. Altogether, we deduce that uncountably many actions from $\Cal F$ are non-orbit equivalent.

\vskip 0.1in
{\it Added in the proof.} Recently, a combination of results and ideas from [E], [GL] and the present paper has been used to show that the above result holds true for any non-amenable group $\Gamma$ ([E]). This development led us to several applications which we present in Section 4.

\vskip 0.2in
 \centerline{\S1. {\bf A separability argument.}}
\vskip 0.1in
\noindent
{\bf 1.1. Conventions.} We start this section by recalling some of the notions that we will use further. 
 For this, fix two m.p actions $\Gamma\curvearrowright^{\sigma}(X,\mu)$ and $\Gamma\curvearrowright^{\alpha}(Z,\nu)$ of a countable group $\Gamma$. 

\vskip 0.03in
$(i)$ The {\bf unitary representation} $\pi_{\sigma}:\Gamma\rightarrow\Cal U((L^2(X,\mu))$ {\bf induced by $\sigma$} is  defined by $\pi_{\sigma}(\gamma)(f)=f\circ\sigma(\gamma^{-1}),$ for all $f\in L^2(X,\mu)$ and  $\gamma\in\Gamma$. We denote by $\pi_{\sigma}^0$ the restriction of $\pi_{\sigma}$ to $L^2(X,\mu)\ominus\Bbb C1$.

$(ii)$ If $Y\subset X$ is a measurable $\sigma(\Gamma)$-invariant set, then we call the action $\Gamma\curvearrowright^{\sigma}(Y,\frac{\mu_{|Y}}{\mu(Y)})$ the {\bf restriction of $\sigma$ to $Y$} and we denote it $\sigma_{|Y}$. In this case, we have that $\pi_{\sigma}=\pi_{{\sigma}_{|Y}}\oplus\pi_{{\sigma}_{|X\setminus Y}}$.

$(iii)$ We say that $\alpha$ is a {\bf quotient} of $\sigma$ if there exists a measurable, measure preserving, onto map $p:X\rightarrow Z$ (called {\bf the quotient map}) such that $p\circ\sigma(\gamma)=\alpha(\gamma)\circ p$, for all $\gamma\in\Gamma$. In this case, we have that $\pi_{\alpha}\subset\pi_{\sigma}$.

$(iv)$ We say that $\alpha$ and $\sigma$ are {\bf conjugate} if there exists a measure space isomorphism $p:X\rightarrow Z$ satisfying the condition in $(iii)$. In this case, we have that $\pi_{\alpha}=\pi_{\sigma}$.

$(v)$ The {\bf diagonal product} of $\alpha$ and $\sigma$ is the action of $\Gamma$ on $(Z,\nu)\times (X,\mu)$ given by $(\alpha\times\sigma)(\gamma)=\alpha(\gamma)\times\sigma(\gamma)$ for all $\gamma\in\Gamma$. 
\vskip 0.1in
\noindent
{\bf 1.2. Relative property (T).}
For an inclusion $\Gamma_0\subset\Gamma$ of countable, discrete  groups  we say that the pair $(\Gamma,\Gamma_0)$ has {\bf relative property} (T) if
for all $\varepsilon>0$, there exists $\delta>0$ and $F\subset \Gamma$ finite such that if $\pi:\Gamma\rightarrow\Cal U(\Cal H)$ is a unitary representation and $\xi\in \Cal H$ is a unit vector satisfying $$||\pi(g)(\xi)-\xi||<\delta,\forall g\in F,$$ then there exists $\xi_0\in\Cal H$ such that $$||\xi_0-\xi||<\varepsilon, \pi(h)(\xi_0)=\xi_0,\forall h\in \Gamma_0.$$ Following Kazhdan-Margulis, the pair (SL$_2(\Bbb Z)\ltimes\Bbb Z^2,\Bbb Z^2)$, where  SL$_2(\Bbb Z)$ acts on $\Bbb Z^2$ by matrix multiplication, has relative property (T) ([K],[M]).  For more examples of pairs of groups with relative property (T), see [Sh1] and [Fe].

\vskip 0.05in
 From now on, we fix a countable group $\Gamma$ which contains a copy  of $\Bbb F_2$. We also fix a free subgroup $\Bbb F_2\subset \Gamma$. 
 Next, if we view $\Bbb F_2$ as a finite index subgroup of SL$_2(\Bbb Z)$.  In particular, we get that the pair $(\Bbb F_2\ltimes\Bbb Z^2,\Bbb Z^2)$ has relative property (T). In fact, for any non-amenable subgroup $\Gamma$ of SL$_2(\Bbb Z)$, the pair $(\Gamma\ltimes\Bbb Z^2,\Bbb Z^2)$ has relative property (T) ([Bu]). 
 Also, we denote by $\alpha$ the action of $\Bbb F_2$ on $\Bbb T^2=\hat{\Bbb Z^2}$ induced by the action of $\Bbb F_2$ on $\Bbb Z^2$. Note that this action preserves the Haar measure $\lambda^2$ of $\Bbb T^2$ and is free and weakly mixing. Finally, we represent the group $\Bbb F_2\ltimes\Bbb Z^2$ as $\{(a,\gamma)|a\in \Bbb Z^2,\gamma\in\Bbb F_2\}$ with the group multiplication given by $(a_1,\gamma_1)\circ (a_2,\gamma_2)=(a_1\gamma_1(a_2),\gamma_1\gamma_2)$.
\vskip 0.05in

\proclaim {1.3. Theorem } Let $\Cal F$ be the class of  free, ergodic, m.p. actions $\Gamma\curvearrowright^{\sigma}(X,\mu)$ on a standard probability space $(X,\mu)$ satisfying the following:

\vskip 0.03in
$(i)$ $\alpha$ is a quotient of $\sigma_{|\Bbb F_2}$, with the quotient map $p:X\rightarrow \Bbb T^2$.

\vskip 0.03in
$(ii)$  $\forall \gamma\in\Gamma\setminus\{e\},$ the set $\{x\in X| p(\gamma x)=p(x)\}$ has zero measure.
\vskip 0.03in
Let $\{\sigma_i\}_{i\in I}\subset\Cal F$ be an uncountable family of mutually orbit equivalent  actions. Then there exists an uncountable set $J\subset I$ with the following property: for every $i,j\in J$, there exist two measurable sets $X_i,X_j\subset X$ of positive measure such that $X_i$ is $\sigma_i(\Bbb F_2)$-invariant, $X_j$ is $\sigma_j(\Bbb F_2)$-invariant and the restriction of ${\sigma_i}_{|\Bbb F_2}$ to $X_i$ is conjugate to the restriction of ${\sigma_j}_{|\Bbb F_2}$ to $X_j$.
\endproclaim

{\it Proof.}  Using the hypothesis we can actually assume that all $\sigma_i$ generate the same measurable equivalence relation $\Cal R\subset X\times X$, i.e. $$\Cal R=\{(x,{\sigma_i}(\gamma)(x))|x\in X,\gamma\in\Gamma\}, \forall i\in I.$$  Following [FM], we endow $\Cal R$ with the measure $\tilde\mu$ given by $$\tilde\mu(A)=\int_{X}|A\cap (\{x\}\times X)| d\mu(x),$$ for every Borel subset $A\subset\Cal R.$

By condition $(i)$, for every $j\in I$, we can find  a quotient map $p_j:(X,\mu)\rightarrow (\Bbb T^2,\lambda^2)$  a such that $p_j\circ{\sigma_j}(\gamma)=\alpha(\gamma)\circ p_j,$ for all $\gamma\in\Bbb F_2$. 
If $a\in \Bbb Z^2$, then we view $a$ as a character on $\Bbb T^2$   and we define $\eta_a^j=a\circ p_j\in L^{\infty}(X,\mu),$ for all $j\in I.$ It is easy to see that for all $a\in \Bbb Z^2,\gamma\in\Bbb F_2$ and  $j\in I$ we have that $\eta_{\gamma(a)}^j=\eta_a^j\circ\sigma_j(\gamma^{-1}).$
 Using this relation it follows that the formula $$\pi_{i,j}(a,\gamma)(f)(x,y)=\eta_a^i(x)\overline{\eta_a^j(y)}f(\sigma_i(\gamma^{-1})(x),\sigma_j(\gamma^{-1})(y)),$$ for all $f\in L^2(\Cal R,\tilde\mu),$  $(x,y)\in\Cal R$ and  $(a,\gamma)\in\Bbb F_2\ltimes \Bbb Z^2$,  defines a unitary representation $\pi_{i,j}:\Bbb F_2\ltimes \Bbb Z^2\rightarrow \Cal U(L^2(\Cal R,\tilde\mu)),$ for every $i,j\in I$.
\vskip 0.05in
 Let $\xi=1_{\Delta}$, where $\Delta=\{(x,x),x\in X\}$, then $\xi\in L^2(\Cal R,\tilde\mu)$ and $||\xi||_{L^2(\Cal R,\tilde \mu)}=1$. For every $i,j\in I$ and all $(a,\gamma)\in\Bbb F_2\ltimes\Bbb Z^2$ we have that 
$$||\pi_{i,j}(a,\gamma)(\xi)-\xi||_{L^2(\Cal R,\tilde\mu)}^2=2-2 \Re <\pi_{i,j}(a,\gamma)(\xi),\xi>_{L^2(\Cal R,\tilde\mu)}=\tag 1$$
$$2-2 \Re \int_{X} \eta_a^i(x)\overline{\eta_a^j(x)}1_{\{{\sigma_i}(\gamma^{-1})= {\sigma_j}(\gamma^{-1})\}}(x)d\mu(x)\leq$$ ( if $||f||_{\infty},||g||_{\infty}\leq 1$, then $\Re\int_{X}(1-fg)\leq ||1-f||_2+||1-g||_2)$ 

$$2||1-1_{\{{\sigma_i}(\gamma^{-1})={\sigma_j}(\gamma^{-1})\}}||_{L^2(X,\mu)}+2||1-\eta_a^i\overline{\eta_a^j}||_{L^2(X,\mu)}=$$ $$2||1_{\{(x,{\sigma_i}(\gamma^{-1})(x))|x\in X\}}-1_{\{(x,{\sigma_j}(\gamma^{-1})(x))|x\in X\}}||_{L^2(\Cal R,\tilde \mu)}+2||\eta_a^i 1_{\Delta}-\eta_a^j 1_{\Delta}||_{L^2(\Cal R,\tilde\mu)}.$$
\vskip 0.05in
Now, since the pair $(\Bbb F_2\ltimes\Bbb Z^2,\Bbb Z^2)$ has relative property (T), we can find  $\delta>0$ and $A\subset \Bbb Z^2,B\subset \Bbb F_2$ finite sets such that if $\pi:\Bbb F_2\ltimes\Bbb Z^2\rightarrow \Cal U(\Cal H)$ is a unitary representation and $\xi\in \Cal H$ is a unit vector which satisfies $$||\pi(a,\gamma)(\xi)-\xi||\leq \delta, \forall (a,\gamma)\in A\times B,$$ then there exists $\xi_0\in\Cal H$ a $\pi(\Bbb Z^2)$-invariant vector such that $||\xi_0-\xi||\leq 1/2$.  

Next, since the Hilbert space $L^2(\Cal R,\tilde\mu)$ is separable and $I$ is uncountable,  we can find $J\subset I$ uncountable such that $$||\eta_a^i\chi_{\Delta}-\eta_a^j\chi_{\Delta}||_{L^2(\Cal R,\tilde\mu)}\leq \delta^2/4,\forall a\in A $$ and $$||\chi_{\{(x,{\sigma_i}(\gamma^{-1})(x))|x\in X\}}-\chi_{\{(x,{\sigma_j}(\gamma^{-1})(x))|x\in X\}}||_{L^2(\Cal R,\tilde \mu)}\leq \delta^2/4,\forall \gamma\in B,$$ for all $i,j\in J$. When combined with inequality (1), this gives that $$||\pi_{i,j}(a,\gamma)(\xi)-\xi||_{L^2(\Cal R,\tilde\mu)}\leq\delta,\forall (a,\gamma)\in A\times B,\forall i,j\in J\tag 2$$ 
Fix $i,j\in J$. Then by using relative property (T) together with (2)  we can find $f\in L^2(\Cal R,\tilde\mu)$ such that $||f-\chi_{\Delta}||_{L^2(\Cal R,\tilde\mu)}\leq 1/2$ and $f$ is $\pi_{i,j}(\Bbb Z^2)$-invariant, i.e. $$ f(x,y)=\eta_a^i(x)\overline{\eta_a^j(y)}f(x,y),\forall a\in \Bbb Z^2\tag 3$$ $\tilde\mu$ almost everywhere $(x,y)\in\Cal R$.
\vskip 0.05in
 Define $$S=\{(x,y)\in\Cal R|\eta_a^i(x)=\eta_a^j(y),\forall a\in \Bbb Z^2\}.$$ Then $S\subset\Cal R$ is measurable and since $f\not=0$, (3) implies that $\tilde{\mu}(S)>0$.
 We claim that for almost every $x\in X$, there is at most one $y\in X$  such that $(x,y)\in S$. If not, then we can find $X_0\subset X$ a set of positive measure and $\gamma\not=\gamma'\in\Gamma$ such that $$(x,\sigma_j(\gamma)(x)), (x,\sigma_j(\gamma')(x))\in S,\forall x\in X_0.$$ Thus, in particular, we get that $$\eta_a^j(\sigma_j(\gamma)(x))=\eta_a^j(\sigma_j(\gamma')(x)),\forall x\in X_0,\forall a\in \Bbb Z^2,$$ or, equivalently, $$a(p_j(\sigma_j(\gamma)(x)))=a(p_j(\sigma_j(\gamma')(x))),\forall a\in \Bbb Z^2,\forall x\in X_0.$$ Since characters separate points, we deduce that $p_j(\sigma_j(\gamma)(x))=p_j(\sigma_j(\gamma')(x)),$ for all $x\in X_0$. However, since $X_0$ is assumed to have positive measure, this contradicts condition $(ii)$, thus proving the claim.
Now, define $X_i$  to be the set of $x\in X$ with the property that there exists a unique $y\in X$ such that $(x,y)\in S$.
The above claim and the fact that $\tilde{\mu}(S)>0$ imply that $\mu(X_i)>0$. 

\vskip 0.05in
 If $(x,y)\in S$, then $\eta_a^i(x)=\eta_a^j(y),$ for all $a\in \Bbb Z^2$, thus $$\eta_{\gamma(a)}^i(x)=\eta_{\gamma(a)}^j(y),\forall a\in \Bbb Z^2,\forall \gamma\in\Bbb F_2.$$ Since $\eta_{\gamma(a)}^i=\eta_a^i\circ{\sigma_i}(\gamma^{-1}),$ for all $a\in \Bbb Z^2$ and $\gamma\in\Bbb F_2,$ we deduce that $$({\sigma_i}(\gamma)(x),{\sigma_j}(\gamma)(y))\in S,\forall \gamma\in\Bbb F_2\tag 4$$

In particular, we get that $X_i$ is $\sigma_i(\Bbb F_2)$-invariant. If we denote $X_j=\{y\in X|\exists x\in X_i, (x,y)\in S\}$, then $X_j$ is a measurable $\sigma_j(\Bbb F_2)-$invariant set.  Define $\phi:X_i\rightarrow X_j$ by $y=\phi(x)$ {\it iff} $(x,y)\in S$.
Then $\phi$ is a measure preserving isomorphism. Indeed, as above, it follows that for almost every $y\in X$, there exists at most one $x\in X$ such that $(x,y)\in S$, hence $\phi$ is an isomorphism. Moreover, since $\phi(x)$ lies in the orbit of $x$ for almost every $x\in X$, we get that $\phi$ is measure preserving.
  
Finally, note that relation (4) implies that ${\sigma_j}(\gamma)(\phi(x))=\phi({\sigma_i}(\gamma)(x))$ almost everywhere $x\in X_i$ and for all $\gamma\in\Bbb F_2$, which gives the desired conjugacy.\hfill $\square$
\vskip 0.05in
Note that up to this point we have no examples of class $\Cal F$ actions. This will be done in the next section by using a co--inducing construction for actions.
\vskip 0.2in 
\centerline {$\S$2. {\bf The co--induced action.}}
\vskip 0.1in
 Let $\Gamma_0\subset \Gamma$ be two countable groups and let $\Gamma_0\curvearrowright^{\alpha}(Y,\nu)$ be a m.p. action. 
Then there is a natural way to construct a m.p. action of $\Gamma$ whose restriction to $\Gamma_0$ admits $\alpha$ as a quotient. 
This construction first appeared in [L] (see also [DGRS] and [G]). We would like to mention that we initially learned of this construction from Section 3.4. in [G]. 
Start by defining $$X=\{f:\Gamma\rightarrow Y|f(\gamma\gamma_0)=\alpha(\gamma_0)(f(\gamma)),\forall \gamma_0\in\Gamma_0,\forall\gamma\in\Gamma\}.$$
and note that $\Gamma$ acts on $X$ by the formula $(\gamma f)(\gamma')=f(\gamma^{-1}\gamma'),$ for all $\gamma$ and $\gamma'\in\Gamma.$\vskip 0.05in
Let $e\in S\subset \Gamma$ be a set such that $\Gamma=\sqcup_{s\in S}s\Gamma_0$.
We observe that $X$ can be identified with $Y^S=\prod_{s\in S}Y$ via $f\rightarrow (f(s))_{s\in S}.$ 
Using this identification we get an action $\tilde\alpha$ (called the {\bf co-induced action}
) of $\Gamma$ on $Y^S$ given by $\tilde\alpha(\gamma)((x_s)_s)=(y_{s'})_{s'},$ where $y_{s'}=\alpha(\gamma_0^{-1})(x_s)$ for the unique $s\in S$ and $\gamma_0\in\Gamma_0$ such that $\gamma^{-1}s'=s\gamma_0$.
Then $\tilde\alpha$ preserves the product measure $\nu_S=\otimes_{s\in S}\nu$ on $Y^S$.

 In the next two lemmas we discuss the freeness and ergodicity of $\tilde\alpha$. Before this, we remark that $p:Y^{S}\rightarrow Y$  given by $p((x_s)_s)=x_e$ is a quotient map and that $p$ realizes $\alpha$ as a quotient of ${\tilde\alpha}_{|\Gamma_0}$.
\vskip 0.05in
\noindent
 \proclaim {2.1. Lemma} Assume that $\alpha$ is a free action and that $(Y,\nu)$ is a non-atomic probability space. Then the set $A_{\gamma}=\{x\in Y^S|p(\gamma x)=p(x)\}$ has zero measure, for all $\gamma\in\Gamma\setminus\{e\}$. In particular, $\tilde\alpha$ is free.
\endproclaim
{\it Proof.} Note that if $\gamma\in\Gamma_0\setminus\{e\}$, then $A_{\gamma}=\{x\in Y^S|\gamma x_e=x_e\}$, hence  the freeness of $\alpha$ implies that $A_{\gamma}$ has measure zero.
 On the other hand, if $\gamma\in\Gamma\setminus\Gamma_0$, let $s\in S\setminus\{e\}$  and $\gamma_0\in\Gamma_0$ such that $\gamma^{-1}=s\gamma_0$. Then $A_{\gamma}=\{x\in Y^S|x_e={\gamma_0}^{-1}x_s\}$, and since $Y$ is non-atomic, we get that $\nu_S(A_{\gamma})=0$.\hfill$\square$ 

\vskip 0.1in
\proclaim{2.2. Lemma} In the above setting, let $\Lambda\subset\Gamma$ be  a subgroup. Then 
\vskip 0.01in
 $(i)$  $\tilde\alpha_{|\Lambda}$ is weakly mixing {\it iff} (if and only if) $\alpha_{|s\Lambda s^{-1}\cap \Gamma_0}$ is weakly mixing for any $s\in\Gamma$ such that $s\Lambda s^{-1}\cap \Gamma_0\subset s\Lambda s^{-1}$ is of finite index.
\vskip 0.01in
$(ii)$ If $|\Gamma/\Gamma_0|=\infty$, then $\tilde\alpha$ is weakly mixing. If $|\Gamma/\Gamma_0|<\infty$, then $\tilde\alpha$ is weakly mixing {\it iff} $\alpha$ is weakly mixing. 
\vskip 0.01in
$(iii)$ $\tilde\alpha_{|\Lambda}$ is  mixing {\it iff} $\alpha_{|s\Lambda s^{-1}\cap \Gamma_0}$ is  mixing for any $s\in\Gamma$.
\vskip 0.01in
$(iv)$ $\tilde\alpha$ is mixing {\it iff} $\alpha$ is mixing.
\endproclaim
{\it Proof.}  $(i)$ Consider the action of $\Gamma$ on $S$ given by $$\gamma\cdot s'=s\Leftrightarrow\gamma s'\in s\Gamma_0.$$

 For every $t\in S$ and $\gamma\in\Lambda$, let $\beta_t$ be the m.p. action of $\Lambda$ on  $\prod_{s\in \Lambda\cdot t}(Y,\nu)_s$  given by 
$\beta_t(\gamma)((x_s)_s)=(y_{s'})_{s'},$ where $y_{s'}=\alpha(\gamma_0^{-1})(x_s)$, for the unique $s\in \Lambda\cdot t$ and $\gamma_0\in\Gamma_0$ such that $\gamma^{-1}s'=s\gamma_0$. Note that if $T\subset S$ is such that $S=\sqcup_{t\in T}\Lambda\cdot t$, then $\tilde\alpha_{|\Lambda}$ is the diagonal product of the actions $\beta_t$ with $t\in T$, i.e. $$\tilde\alpha_{|\Lambda}=\times_{t\in T}\beta_t.$$
\vskip 0.05in
 
{\it Claim 1.} If $t\in T$ and $\Lambda\cdot t$ is infinite, then $\beta_t$ is weakly mixing. 
\vskip 0.05in
{\it Proof of Claim 1.}
To prove that $\beta_t$ is weakly mixing  we need to show that if $\xi_1,..,\xi_n\in L^2(\prod_{s\in \Lambda\cdot t}(Y,\nu)_s)$ are functions of zero integral, then for every $\varepsilon>0$ we can find  $\gamma\in\Lambda$ such that $ |<\beta_t(\gamma)(\xi_i),\xi_j>|\leq\varepsilon$ for all $i,j$. Note that in order to prove this condition, we can assume that there exists a finite set $F\subset\Lambda\cdot t$ such that $\xi_i\in L^2(\prod_{s\in F}(Y,\nu)_s)$ for all $i=1,..,n.$

Now, since $\Lambda\cdot t$ is infinite, we can find $\gamma\in\Lambda$ such that $\gamma F\cap F=\emptyset$. This implies that $\beta_t(\gamma)(\xi_i)$ and $\xi_j$ are independent for all $i,j$. Thus, $<\beta_t(\gamma)(\xi_i),\xi_j>=0$ for all $i,j$, hence $\beta_t$ is weakly mixing.

\vskip 0.05in
Using Claim 1 we get that $\tilde\alpha_{|\Lambda}$ is weakly mixing iff $\beta_t$ is weakly mixing for every $t\in T$ such that $\Lambda\cdot t$ is finite.
Let $t\in T$  such that $\Lambda\cdot t$ is finite. Then$$\Lambda_t=\{\gamma\in\Lambda|\gamma\cdot t'=t',\forall t'\in\Lambda\cdot t\}$$ is a finite index subgroup of $\Lambda$.  Thus, $\beta_t$ is weakly mixing iff ${\beta_t}_{|\Lambda_t}$ is weakly mixing.
Since $$\beta_t(\gamma)=\times_{s\in \Lambda\cdot t}\alpha(s^{-1}\gamma s),\forall\gamma\in\Lambda_t,$$ we further deduce that ${\beta_t}$ is weakly mixing iff $\alpha_{|s^{-1}\Lambda_t s}$ is weakly mixing for every $s\in\Lambda\cdot t$.  Next, note that the inclusions $$s^{-1}\Lambda_ts\subset s^{-1}\Lambda s\cap\Gamma_0\subset s^{-1}\Lambda s$$ are of finite index for every $s\in\Lambda\cdot t$. This implies that $\beta_t$ is  weakly mixing iff  $\alpha_{|s^{-1}\Lambda s\cap \Gamma_0}$ is weakly mixing for every $s\in\Lambda\cdot t$. Altogether, we get that $\tilde\alpha_{|\Lambda}$ is weakly mixing iff $\alpha_{|s^{-1}\Lambda s\cap \Gamma_0}$ is weakly mixing for all $s\in\Gamma$ such that $\Lambda\cdot s$ is finite.
\vskip 0.02in
Since $\Lambda\cdot s$ is finite iff $\Lambda\cap s\Gamma_0 s^{-1}=\{\gamma\in\Lambda|\gamma\cdot s=s\}\subset\Lambda$ is of finite index, we get the conclusion.
\vskip 0.02in
$(iii)$ Assume that  $\alpha_{|s\Lambda s^{-1}\cap \Gamma_0}$ is  mixing for any $s\in\Gamma$.
To prove that ${\tilde\alpha}_{|\Lambda}$ is mixing it suffices to show the following:
\vskip 0.05in
{\it Claim 2.}
For any  finitely supported vectors $f=\otimes_{s\in A}f_s, g=\otimes_{s\in B}g_s\in L^{\infty}(Y^S,\nu_S),$ where $A,B\subset S$ are finite and $f_s,g_t\in L^{\infty}(Y,\nu)$ have zero integral, for all $s\in A$ and $t\in B$, we have that $$\lim_{\Lambda\ni\gamma\rightarrow\infty}<\tilde\alpha(\gamma)(f),g>=0.$$
\vskip 0.05in
{\it Proof of Claim 2.}
 Note that the induced action $\tilde\alpha:\Gamma\rightarrow\text{Aut}(L^{\infty}(Y^S))$ is given by $\tilde\alpha(\gamma)(\otimes_s f_s)=\otimes_{s'}g_{s'},$ where $g_{s'}=\alpha(\gamma_0)(f_s)$ for the unique $s\in S$ and $\gamma_0\in\Gamma_0$ such that $\gamma s=s'\gamma_0$. Using this we get that if $\gamma\in\Gamma$, then $<\tilde\alpha(\gamma)(f),g>=0$ unless $|A|=|B|$ and there exists a bijection $\pi:A\rightarrow B$ such that ${\pi(s)}^{-1}\gamma s\in\Gamma_0$ for all $s\in A$. In the latter case, we have that $$<\tilde\alpha(\gamma)(f),g>=\prod_{s\in A}<\alpha({\pi(s)}^{-1}\gamma s))(f_s),g_{\pi(s)}>.$$ For a bijection $\pi:A\rightarrow B$, let  $\Lambda_{\pi}=\{\gamma\in\Lambda|{\pi(s)}^{-1}\gamma s\in\Gamma_0, \forall s\in A\}$. Then, proving the claim is equivalent to proving that $$\lim_{\Lambda_{\pi}\ni\gamma\rightarrow\infty}<\tilde\alpha(\gamma)(f),g>=0,$$ for all bijections $\pi:A\rightarrow B$. Fix a bijection $\pi:A\rightarrow B$, $\lambda\in\Lambda_{\pi}$ and $s\in A$. Then for all $\gamma\in\Lambda_{\pi}$ we have that  $s^{-1}(\lambda^{-1}\gamma)s\in s^{-1}\Lambda s\cap\Gamma_0$ and that $$<\alpha({\pi(s)}^{-1}\gamma s))(f_s),g_{\pi(s)}>=<\alpha(({\pi(s)}^{-1}\lambda s) (s^{-1}(\lambda^{-1}\gamma)s))(f_s),g_{\pi(s)}>=$$ $$<\alpha(s^{-1}(\lambda^{-1}\gamma)s)(f_s), \alpha(s^{-1}\lambda^{-1}\pi(s))(f_{\pi(s)})>$$
Now, if we let $\Lambda_{\pi}\ni\gamma\rightarrow\infty$, then $ s^{-1}\Lambda s\cap\Gamma_0\ni s^{-1}(\lambda^{-1}\gamma)s\rightarrow\infty$.
Since $\alpha_{|s^{-1}\Lambda s\cap\Gamma_0 }$ is mixing by our assumption, we get that $$\lim_{\Lambda_{\pi}\ni\gamma\rightarrow\infty}<\alpha({\pi(s)}^{-1}\gamma s))(f_s),g_{\pi(s)}>=0,$$ which ends the proof of the claim.

The other implication follows easily and we omit its proof. 
Finally, $(ii)$ and $(iv)$ follow by applying $(i)$ and $(iii)$ to $\Lambda=\Gamma$.
\hfill $\square$
\vskip 0.05in 
For the next result, we assume the notations and assumptions of Section 1. Thus, $\Gamma$ is a countable group which contains $\Gamma_0=\Bbb F_2$ and $\alpha$ denotes the action $\Bbb F_2\curvearrowright \Bbb T^2$.   
\proclaim {2.3. Corollary}  $\tilde\alpha$ is weakly mixing and belongs to $\Cal F$.
Moreover, for any ergodic action $\rho$ of $\Gamma$, the diagonal product action $\tilde\alpha\times\rho$ also belongs to $\Cal F$. 
\endproclaim
{\it Proof.} Since $\alpha$ is weakly mixing
(see for example [P2]), Lemma 2.2.$(ii)$ implies that $\tilde\alpha$ is weakly mixing. When combined with Lemma 2.1. this gives that $\tilde\alpha\in\Cal F$.
The second assertion follows easily since $\tilde\alpha$ is weakly mixing and thus the diagonal product with any ergodic action is still ergodic.\hfill $\square$
\vskip 0.2in
\centerline{\S3. {\bf Proof of the Main Result.}}
\vskip 0.1in
 Let $\Gamma$ be a countable group containing $\Bbb F_2$. 
Let $\{\pi_i:\Bbb F_2\rightarrow\Cal U(\Cal H_i)\}_{i\in I}$ be an uncountable family of mutually non-equivalent, irreducible, {\bf c}$_0$-representations of $\Bbb F_2$, i.e. such that $\lim_{g\rightarrow\infty}<\pi_i(g)\xi,\eta>=0$ for all $\xi,\eta\in\Cal H_i$ and $i\in I$ ([Sz]). 
\vskip 0.05in
 {\it  Claim 1.} For every $i\in I$, there exists a free, mixing, m.p. action $\Gamma\curvearrowright^{\rho_i}(X_i,\mu_i)$ such that  $$\pi_i\subset {\pi_{\rho_i}^0}_{|\Bbb F_2}.$$ 
\vskip 0.05in
{\it Proof of Claim 1.}
For every $i\in I$, let $\tilde{\pi}_i:\Gamma\rightarrow\Cal U(\tilde \Cal H_i)$ be the induced representation. 
 Then we can find a free, m.p., Gaussian action $\Gamma\curvearrowright^{\rho_i}(X_i,\mu_i)$ (see for example [Fu2],[Ke]) such that 
$$\tilde\pi_i\subset\pi_{\rho_i}^0\subset \oplus_{n\geq 1}{\tilde\pi_i}^{\otimes_n}.$$
Now, since $\pi_i$ is {\bf c}$_0$, we get that $\tilde\pi_i$ is also {\bf c}$_0$, thus $\rho_i$ is a mixing action. Also, since $\pi_i\subset \tilde{\pi_i}_{|\Bbb F_2}$ we get the second assertion.
\vskip 0.1in

Next, fix an ergodic, measure preserving  action $\Gamma\curvearrowright^{\beta} (Y,\nu)$. For every $i\in I$, consider the diagonal product action $\sigma_i=\tilde\alpha\times\beta\times\rho_i$ of $\Gamma$ on $$(Z_i,\eta_i):=\prod_{s\in \Gamma/\Bbb F_2}(\Bbb T^2,\lambda^2)_s\times (Y,\nu)\times (X_i,\mu_i),$$ where $\tilde\alpha$ denotes the action $\Gamma\curvearrowright \prod_{s\in \Gamma/\Bbb F_2}(\Bbb T^2,\lambda^2)_s$ obtained by co-inducing $\alpha$.
Since $\tilde\alpha$ is free, weakly mixing, $\beta$ is ergodic and $\rho_i$ is mixing, we deduce that $\sigma_i$ is a free, ergodic action, for all $i\in I$. 
 \vskip 0.05in
{\it Claim 2}. Let $i\in I$ and let $Z_i'\subset Z_i$ be a $\sigma_i(\Bbb F_2)-$invariant set. Then the representation induced by the restriction of   ${\sigma_i}_{|\Bbb F_2}$ to $Z_i'$ contains $\pi_i$.
\vskip 0.05in
{\it Proof of Claim 2}. Since $\rho_i$ is mixing, we derive that ${\rho_i}_{|\Bbb F_2}$ is weakly mixing. Thus, since $Z_i'$ is $\sigma_i(\Bbb F_2)$-invariant, we get that $Z_i'= B\times X_i,$ for some measurable  set $B\subset \prod_{s\in \Gamma/\Bbb F_2}\Bbb T^2\times Y$. This implies that the restriction of ${\sigma_i}_{|\Bbb F_2}$ to $Z_i'$ admits ${\rho_i}_{|\Bbb F_2}$ as a quotient. Thus, the representation induced by the restriction of   ${\sigma_i}_{|\Bbb F_2}$ to $Z_i'$ contains ${\pi_{\rho_i}^0}_{|\Bbb F_2}$.
Since, by Claim 1, the latter contains $\pi_i$, we are done.

\vskip 0.1in
{\it Claim 3.} For every $i\in I$, the set $I_i$ of $j\in I$ such that a restriction of ${\sigma_j}_{|\Bbb F_2}$ is conjugate to  a restriction of ${\sigma_{i}}_{|\Bbb F_2}$ is countable.
\vskip 0.05in
{\it Proof of Claim 3.} Let $\pi$ be the unitary representation of $\Bbb F_2$ induced by ${\sigma_{i}}_{|\Bbb F_2}$. If $j\in I_i$, then $\pi$ contains the representation induced by a restriction of ${\sigma_j}_{|\Bbb F_2}$. Now, by  Claim
2, the latter contains $\pi_j$ as a subrepresentation. Combining these two inclusions, we get that $\pi_j\subset \pi,$ for all $j\in I_i.$
Since a separable unitary representation can only have countably many non-equivalent irreducible subrepresentations and since the $\pi_j's$ are irreducible and mutually non-equivalent, it follows that $I_i$ is countable.  
\vskip 0.1in
 
{\it Claim 4.} Uncountably many of the actions $\{\sigma_i\}_{i\in I}$ are mutually non-OE. 
\vskip 0.05in
{\it Proof of Claim 4.} If we assume the contrary, then we can find an uncountable set $J\subset I$ such that the actions $\{\sigma_j\}_{j\in J}$ are mutually orbit equivalent. Now, since $\beta\times\rho_i$ is ergodic, Corollary 2.3. implies that $\sigma_i\in\Cal F,$ for all $i\in I.$   Thus, by applying Theorem 1.3. to the family of actions $\{\sigma_j\}_{j\in J}\subset \Cal F$, we can find an uncountable subset $K\subset J$ such that for all $k,l\in K$, a restriction of ${\sigma_k}_{|\Bbb F_2}$ is conjugate to a restriction of ${\sigma_l}_{|\Bbb F_2}$. This, however, implies that $I_k$ is uncountable, for every $k\in K$, in contradiction with Claim 3. 
\vskip 0.05in
Now, if we take $\beta$ to be any ergodic action of $\Gamma$ (e.g. the trivial action), then Claim 4 gives uncountably many non-OE actions of $\Gamma$. 
\vskip 0.05in

For the  moreover assertion, we first recall the definition of measure equivalence (see [Fu1]).
To this end, let $\Gamma\curvearrowright^{\beta}(Y,\mu)$ be a free,  ergodic, m.p. action   and let $t>0$. Let $n>t$ be a natural number and set $Y^n=Y\times\{1,..,n\}$ endowed with the natural measure.  Next, let $Y^t\subset Y^n$ be a measurable set of measure $t$ and define $\Cal R_{\beta}^t$ be the  equivalence relation on $Y^t$ given by: $(x,i)\sim(y,j)$ {\it iff} there exists $\gamma\in\Gamma$ such that $y=\beta(\gamma)(x)$. Note that the isomorphism class of $\Cal R_{\beta}^t$ depends on $t$ and not on the particular choice of $Y^t$ (since $\beta$ is ergodic). If $t=1,$ then we use the notation $\Cal R_{\beta}$.
Two groups $\Gamma$ and $\Lambda$ are {\bf measure equivalent (ME)} if we can find  an action $\beta$  as above, $t>0$ and a free, ergodic, m.p. action $\delta$ of $\Lambda$ on $Y^t$ such that $$\Cal R_{\beta}^t=\Cal R_{\delta}.$$
\vskip 0.05in
Now, let  $\Gamma\curvearrowright^{\theta} (S,m)$ be a weakly mixing, m.p. action. Then the diagonal product action $\theta\times\beta$ is ergodic and the equivalence relation $\Cal R_{\theta\times\beta}^t$ can be realized as the equivalence relation on $S\times Y^t$ given by: $(s,((x,i))\sim (s',(y,j))$ {\it iff} there exists $\gamma\in\Gamma$ such that $s'=\theta(\gamma)(s)$ and $y=\beta(\gamma)(x)$.
Next, we note the following claim  due to Gaboriau ([G]):
\vskip 0.05in
{\it Claim 5.} In the context from above, there exists a free, ergodic, m.p. action $\tau$ of $\Lambda$ on $S\times Y^t$ such that $$\Cal R_{\theta\times\beta}^t=\Cal R_{\tau}.$$ 
\vskip 0.05in
{\it Proof of Claim 5.} Let $\lambda\in\Lambda$, then for almost every $(x,i)\in Y^t$ we can find a unique (by the freeness of $\beta$) $\gamma=w(\lambda,(x,i))\in\Gamma$ such that $\delta(\lambda)(x,i)=(\beta(\gamma)(x),i)$. Then $w:\Lambda\times Y^t\rightarrow\Gamma$ gives a cocycle for $\delta$.
 This implies that the formula $$\tau(\lambda) (s,(x,i))=(\theta(w(\lambda,(x,i)))(s),\delta(\lambda)(x,i))$$ for all $\lambda\in\Lambda$, $s\in S,$ $(x,i)\in Y^t$ defines a m.p. $\Lambda$-action on $S\times Y^t$.
 Moreover, it is clear that $\Cal R_{\theta\times\beta}^t=\Cal R_{\tau},$ hence, since $\theta\times\beta$ is ergodic, we get that $\tau$ is also ergodic. Also, since $\tau$ admits $\delta$ as a quotient and since $\delta$ is free, we deduce that $\tau$ is  free. 
\vskip 0.05in
Finally, let $\Gamma$ be a group containing $\Bbb F_2$ and let $\Lambda$ a group {\bf ME} to $\Gamma$. Let $\beta$ (resp. $\delta$) be a free, ergodic, m.p. action of $\Gamma$ (resp. of $\Lambda$) such that $\Cal R_{\beta}^t=\Cal R_{\delta},$ for some $t>0$.
For all $i\in I$, denote $\theta_i=\tilde\alpha\times\rho_i$ and observe that $\sigma_i=\theta_i\times\beta$.
  By applying Claim 5, it follows that for every $i\in I$ there exists a free, ergodic, m.p. action $\tau_i$ of $\Lambda$ such that $\Cal R_{\sigma_i}^t=\Cal R_{\tau_i}.$ 
Recall that two actions are orbit equivalent {\it iff} they generate isomorphic equivalence relations. Thus, Claim 4 implies that uncountably many of the actions $\{\tau_i\}_{i\in I}$ are mutually non-orbit equivalent. \hfill $\square$ 
\vskip 0.2in

\centerline{\S4. {\bf Applications to von Neumann algebras.}}
\vskip 0.1in
\noindent
{\bf (I)}
After the first draft of this paper was posted on the arxiv (January 2007), there have been two important developments, in [GL] and [E].
 To briefly present these results, recall first that, in general, a non-amenable group $\Gamma$ need not contain a copy of $\Bbb F_2$ ([O]). Nevertheless, D. Gaboriau and R. Lyons proved in [GL] that {\bf any} non-amenable group $\Gamma$ admits $\Bbb F_2$  as a {\it measurable subgroup}:

\proclaim {4.1. Theorem [GL]} Let $\Gamma$ be a countable non-amenable group. Then there exist free, ergodic, m.p. actions $\Gamma\curvearrowright (Z,\eta)$ and $\Bbb F_2\curvearrowright (Z,\eta)$ such that $\Bbb F_2z\subset\Gamma z$, a.e. $z\in Z$.
\endproclaim

This result opened up the possibility that the condition {\it $\Gamma$  contains a copy of $\Bbb F_2$} in the statement of our main theorem could be replaced by the more general, natural condition {\it $\Gamma$ is non-amenable}. To do this, by analogy with the proof of our main result, 
 a co-inducing construction in a group/measurable subgroup situation rather than in a group/subgroup one, was needed.

Recently, I. Epstein obtained such a construction in [E] (see Lemma 4.2.).  Using this construction, she was able to push our arguments in the case that $\Gamma$ is an arbitrary non-amenable group and to show that indeed any such $\Gamma$ admits uncountably many non-OE actions ([E]).

\proclaim{ 4.2. Lemma [E]} Let $\Gamma_0,\Gamma$ be two countable groups and assume that there exist free, ergodic, m.p. actions $\Gamma\curvearrowright (Z,\eta)$ and $\Gamma_0\curvearrowright (Z,\eta)$ such that $\Gamma_0 z\subset\Gamma z$, a.e. $z\in Z$. Let $\Gamma_0\curvearrowright^{\alpha}(Y,\nu)$ be  a free, ergodic, m.p. action. Then there exist a probability space $(X,\mu)$, a quotient map $p:X\rightarrow Y$ and free, ergodic, m.p. actions $\Gamma_0\curvearrowright^{\beta}(X,\mu)$, $\Gamma\curvearrowright^{\tilde\alpha}(X,\mu)$ such that 
\vskip 0.03in
(i)  $\alpha$ as a quotient of ${\beta}$ with $p$ as the quotient map.
\vskip 0.03in
(ii) $\forall\gamma\in\Gamma\setminus\{e\}$, the set $\{x\in X|p(\gamma x)=p(x)\}$ has zero measure.
\vskip 0.03in
(iii) $\Gamma_0 x\subset\Gamma x$, a.e. $x\in X$.
\endproclaim

Below, we obtain consequences of 4.1. and 4.2. in the theory of von Neumann algebras.  We note that in the first draft of this paper, we obtained these corollaries under the additional assumption that $\Gamma$ contains a copy of $\Bbb F_2$.

\vskip 0.1in
\noindent
 {\bf (II)} We begin by observing that one can characterize the non-amenability of a group $\Gamma$ in terms of Popa's notion of relative property (T) for von Neumann algebras. 
For this, let $M$ be a finite  von Neumann algebra with a faithful, normal trace $\tau$ and let $B\subset M$ be a von Neumann subalgebra. The inclusion $(B\subset M)$ is {\bf rigid} (or has {\bf relative property (T)}) if whenever $\phi_n:M\rightarrow M$ is a sequence of unital, tracial, completely positive (c.p.) maps such that $\phi_n\rightarrow$id$_{M}$ in the pointwise $||.||_2$-topology, we must have that $\lim_{n\rightarrow\infty}\sup_{x\in B,||x||\leq 1}||\phi_n(x)-x||_2=0$ ([P1]). In the case $(B\subset M)=(L(\Gamma_0)\subset L(\Gamma))$, for two countable groups $\Gamma_0\subset \Gamma$, the inclusion $(B\subset M)$ is rigid {\it iff} the pair $(\Gamma,\Gamma_0)$ has relative property (T) ([P1]).

Also, recall that the {\bf group measure space} construction  associates that to every free, ergodic, m.p. action $\Gamma\curvearrowright^{\sigma} (X,\mu)$  a II$_1$ factor, $L^{\infty}(X,\mu)\rtimes_{\sigma}\Gamma$, together with a Cartan subalgebra, $L^{\infty}(X,\mu)$ ([MvN]). 
In [P1], Popa asked to characterize the countable groups $\Gamma$ which admit a {\bf rigid} action $\sigma$, i.e. such that the inclusion $L^{\infty}(X,\mu)\subset L^{\infty}(X,\mu)\rtimes_{\sigma}\Gamma$ is rigid. The following result is motivated by this question.

\proclaim {4.3. Theorem} A countable group $\Gamma$ is non-amenable if and only if there exists a free, ergodic, m.p. action $\Gamma\curvearrowright (X,\mu)$ and a diffuse von Neumann subalgebra $Q\subset L^{\infty}(X,\mu)$ such that
\vskip 0.03in
(i) $Q'\cap L^{\infty}(X,\mu)\rtimes\Gamma=L^{\infty}(X,\mu)$ and
\vskip 0.03in 
(ii) the inclusion $Q\subset L^{\infty}(X,\mu)\rtimes\Gamma$ is rigid.
\endproclaim
{\it Proof.} If $\Gamma$ is amenable, then  $L^{\infty}(X,\mu)\rtimes\Gamma$ is isomorphic to the hyperfinite II$_1$ factor, $R$, for any free, ergodic, m.p. action $\Gamma\curvearrowright (X,\mu)$ ([C1],[OW]). Since $R$ does not contain any diffuse von Neumann subalgebra with the relative property (T), we get the "if" part of the conclusion.

For the converse, let $\Gamma$ be a non-amenable group. As before, denote by $\alpha$ the action $\Bbb F_2\curvearrowright (\Bbb T^2,\lambda^2)$. Then, by combining 4.1. and 4.2. we can find a probability space $(X,\mu)$, a quotient map $p:X\rightarrow Y=\Bbb T^2$ and two free, ergodic, m.p. actions $\Bbb F_2\curvearrowright^{\beta}(X,\mu)$, $\Gamma\curvearrowright^{\tilde\alpha} (X,\mu)$ which satisfy conditions $(i)-(iii)$ in 4.2.

Denote by $\theta: L^{\infty}(\Bbb T^2,\lambda^2)\hookrightarrow L^{\infty}(X,\mu)$ the embedding given by $\theta(f)=f\circ p$, for all $f\in L^{\infty}(\Bbb T^2,\lambda^2)$, and let
$Q=\theta(L^{\infty}(\Bbb T^2,\lambda^2))$. We claim that $\tilde\alpha$ and $Q$ verify the conclusion. For this, denote by $\{u_{\gamma}\}_{\gamma\in\Gamma}$ the canonical unitaries implementing the action of $\Gamma$ on $L^{\infty}(X,\mu)$. Then it is easy to see that $Q'\cap L^{\infty}(X,\mu)\rtimes_{\tilde\alpha}\Gamma$ is generated by $L^{\infty}(X,\mu)$ and $\{1_{A_{\gamma}}u_{\gamma}|\gamma\in\Gamma\}$, where $A_{\gamma}=\{x\in X|p(\gamma^{-1} x)=p(x)\}$, for all $\gamma\in\Gamma$. As $\mu(A_{\gamma})=0$, for all $\gamma\in\Gamma\setminus\{e\}$, we deduce that $Q'\cap L^{\infty}(X,\mu)\rtimes_{\tilde\alpha}\Gamma=L^{\infty}(X,\mu)$.

Next, since $p$ realizes $\alpha$ as a quotient of $\beta$, we get that $\theta$ extends to an embedding $$\theta:L^{\infty}(\Bbb T^2,\lambda^2)\rtimes_{\alpha}\Bbb F_2 \hookrightarrow L^{\infty}(X,\mu)\rtimes_{\beta}\Bbb F_2.$$ Now, by [P1], the inclusion $L^{\infty}(\Bbb T^2,\lambda^2)\subset L^{\infty}(\Bbb T^2,\lambda^2)\rtimes_{\alpha}\Bbb F_2$ is rigid, hence the inclusion $Q\subset L^{\infty}(X,\mu)\rtimes_{\beta}\Bbb F_2$ is rigid.
Finally, since $\Bbb F_2x\subset\Gamma x$, a.e. $x\in X$, we have that $ L^{\infty}(X,\mu)\rtimes_{\beta}\Bbb F_2\subset L^{\infty}(X,\mu)\rtimes_{\tilde\alpha}\Gamma$ and we deduce that the inclusion $Q\subset  L^{\infty}(X,\mu)\rtimes_{\tilde\alpha}\Gamma$ is rigid ([P1]).\hfill$\square$
\vskip 0.1in
\noindent
{\bf 4.4. Remark.}  Theorem 4.3. implies that every countable non-amenable group $\Gamma$ admits an {\bf almost rigid} action. To make this precise, let $\sigma$ be the action given by 4.3. and let $\{p_n\}_{n\geq 1}$ be a sequence of projections which generate $L^{\infty}(X,\mu)$. For every $n$, define $Q_n=(Q\vee\{p_1,..,p_n\})''$. Then the inclusion $Q_n\subset L^{\infty}(X,\mu)\rtimes_{\sigma}\Gamma$ is rigid and $Q_n'\cap L^{\infty}(X,\mu)\rtimes_{\sigma}\Gamma= L^{\infty}(X,\mu),$ for all $n$. Moreover, we have that $\overline{\cup_{n\geq 1} Q_n}^{w}=L^{\infty}(X,\mu)$.
\vskip 0.1in

 Next, we denote by $A=\bigoplus_{1}^{\infty}\Bbb Z$ and we note that if $\Gamma$ contains s copy $\Bbb F_2$ then the action $\sigma$  from 4.3. can be taken to come from an action of $\Gamma$ by automorphisms on $A$.
\proclaim {4.5. Proposition} Let $\Gamma$ be a countable group which contains $\Bbb F_2$. Then there exists a homomorphism $\rho: \Gamma\rightarrow \text{Aut}(A)$  and an infinite subgroup $B\subset A$ such that the pair $(\Gamma\ltimes_{\rho}A,B)$ has relative property (T) and that the set $\{\gamma(b)b^{-1}|b\in B\}$ is infinite, for all $\gamma\in\Gamma\setminus\{e\}$.
\endproclaim
{\it Proof.} First, remark that $\Bbb F_2$ contains a copy of itself which has infinite index. Indeed, if $\Bbb F_2=<a,b>$,  then the subgroup generated by $a$ and $bab^{-1}$ has infinite index and is isomorphic to $\Bbb F_2$. Thus, we can assume that $\Bbb F_2$ has infinite index in $\Gamma$. 

 Next, let $e\in S\subset \Gamma$ be a set such that $\Gamma=\sqcup_{s\in S}s\Bbb F_2$ and identify $A$ with $\bigoplus_{s\in S}\Bbb Z^2$.
Then the co-induced construction from Section 2 shows the action $\alpha:\Bbb F_2\rightarrow$ Aut$(\Bbb Z^2)$ co-induces to an action $\rho:\Gamma\rightarrow$ Aut$(A)$. Moreover, we have that $\rho(\Bbb F_2)$ invaries $B=(\Bbb Z^2)_e$ and that the inclusions of groups $(B\subset\Bbb F_2\ltimes_{\rho_{|\Bbb F_2}}B)$ and $(\Bbb Z^2\subset\Bbb F_2\ltimes_{\alpha}\Bbb Z^2)$ are isomorphic. Since the latter inclusion has relative property (T), we deduce that the pair $(\Gamma\ltimes_{\rho}A,B)$ has relative property (T). The second assertion is easy and we leave it to the reader.
\hfill$\square$
\vskip 0.05in
 In the context of 4.5., it now follows that the induced action $\Gamma\curvearrowright^{\sigma}(\hat{A},\mu)$ verifies 4.3, where $\mu$ is the Haar measure on the dual of $A$. We note that we do not know whether the converse of Proposition 4.5. is true, i.e. if any countable group $\Gamma$ which has an action on $A$ with the above properties must necessarily contain $\Bbb F_2$.

\vskip 0.1in

The class of {\bf  $\Cal H\Cal T$ factors} has been introduced by Popa, who used it  to provide the first examples of II$_1$ factors with trivial fundamental group ([P1]).  
A II$_1$ factor $M$ is in the $\Cal H\Cal T$ class if it has a  Cartan subalgebra $A$ (called an HT Cartan subalgebra)  such that: 
\vskip 0.02in
$(i)$ $M$  has the property H relative to $A$ and
\vskip 0.02in
$(ii)$ there exists a von Neumann subalgebra $B\subset A$ such that $B'\cap M\subset A$ and the inclusion $B\subset M$ is rigid.
\vskip 0.05in

 In [P1], Popa raised the question of characterizing HT groups, i.e. groups which admit a free, ergodic, m.p. action $\Gamma\curvearrowright^{\sigma} (X,\mu)$ such that the corresponding Cartan subalgebra inclusion $(A\subset M)=(L^{\infty}(X,\mu)\subset L^{\infty}(X,\mu)\rtimes_{\sigma}\Gamma)$ is HT. 
Since in this case, $M$ has property H relative to $A$ if and only if $\Gamma$ has Haagerup's property ([P1]), Theorem 4.3. implies the following:

\proclaim {4.6. Corollary} A countable group $\Gamma$ is HT if and only if is non-amenable and has Haagerup's  property.
\endproclaim

\vskip 0.05in
\noindent
{\bf (III)} Recall that two actions $\Gamma\curvearrowright (X,\mu)$ and $\Lambda\curvearrowright (Y,\nu)$ are called {\bf von Neumann equivalent (vNE)} if the associated II$_1$ factors are isomorphic, i.e. if $L^{\infty}(X,\mu)\rtimes\Gamma\simeq L^{\infty}(Y,\nu)\rtimes\Lambda$ ([P3]). 
Next, we show that  any non-amenable group admits uncountably many non-von Neumann equivalent actions.
Since von Neumann equivalence of actions is weaker than orbit equivalence ([FM]) this result generalizes Theorem 1.3. as well as the main result of [E].
For $\Gamma=\Bbb F_n, n\geq 2$,
 this result has been first obtained in [GP].

\proclaim{4.7. Theorem} Any countable, non-amenable group $\Gamma$ admits uncountably many non-vNE free, ergodic m.p. actions. 
\endproclaim
{\it Proof.} Let $\tilde{\Cal F}$ be  the class  of free, ergodic, m.p. actions $\Gamma\curvearrowright^{\sigma}(X,\mu)$ such that there exists a free, ergodic, m.p. action $\Bbb F_2\curvearrowright^{\beta}(X,\mu)$ satisfying the following
\vskip 0.03in
$(i)$ $\alpha$ is a quotient of $\beta$, with the quotient map $p:X\rightarrow \Bbb T^2$,
\vskip 0.03in
$(ii)$ $\forall\gamma\in\Gamma\setminus\{e\}$, the set $\{x\in X|p(\gamma x)=p(x)\}$ has zero measure and
\vskip 0.03in
$(iii)$ $\Bbb F_2 x\subset\Gamma x$, a.e. $x\in X$.
\vskip 0.03in 
 It is then proven in [E], by using Theorem 1.3., that there exists an uncountable family of actions $\Gamma\curvearrowright^{\sigma_i}(X_i,\mu_i)$ ($i\in I$) from $\Cal {\tilde F}$ which are mutually non-orbit equivalent.  
For every $i\in I$,  denote $M_i=L^{\infty}(X_i,\mu_i)\rtimes_{\sigma_i}\Gamma$ and $A_i=L^{\infty}(X_i,\mu_i)$. 

\vskip 0.05in
{\it Claim.} For every $i_0\in I$, the set $J=\{i\in I|M_i\simeq M_{i_0}\}$ is countable.
\vskip 0.05in
Note that since $I$ is uncountable, this claim implies that uncountably many of the actions $\{\sigma_i\}_{i\in I}$ are non-von Neumann equivalent.
\vskip 0.05in
{\it Proof of the claim.} Start by denoting $Q=L^{\infty}(\Bbb T^2,\lambda^2)$ and $N=L^{\infty}(\Bbb T^2,\lambda^2)\rtimes_{\alpha}\Bbb F_2$. 
Since $\sigma_i\in\tilde{\Cal F}$, the proof of Theorem 4.3. shows that there exists an embedding  of $N$ into $M_i$ such that under this embedding $Q\subset A_i$ and $Q'\cap M_i=A_i$.
Also, since the inclusion $Q\subset N$ is rigid ([P1]), we can find $F\subset N$ finite and $\delta>0$ such that if a unital, tracial, c.p. map $\phi:N\rightarrow N$ satisfies $||\phi(x)-x||_2\leq \delta$, for all $x\in F$, then $$||\phi(b)-b||_2\leq 1/4,\forall
 b\in (Q)_1.$$

 To prove the claim, assume by contradiction that $J$ is uncountable. For every $i\in J$, let $\theta_i:M_i\rightarrow M_{i_0}$ be an isomorphism and consider the set $\{\theta_i(x)|x\in F\}\subset L^2(M_{i_0})^{\oplus_{|F|}}$. 
    Since $L^2(M_{i_0})$ is a separable Hilbert space and since $J$ is uncountable, we can find $i\not=j\in J$ such that $$||\theta_i(x)-\theta_j(x)||_2\leq \delta,\forall x\in F$$ Thus, the isomorphism $\theta={\theta_j}^{-1}\circ \theta_i:M_i\rightarrow M_j$ satisfies $||\theta(x)-x||_2\leq \delta$, for all $x\in F$. 

Further, if we let $\phi=(E_{N}\circ\theta)_{|N}:N\rightarrow N$ (where $E_N:M_j\rightarrow N$ is the conditional expectation onto $N$), then $\phi$ is a unital, tracial, c.p. map and $$||\phi(x)-x||_2=||E_N(\theta(x))-x||_2=||E_N(\theta(x)-x)||_2\leq \delta,\forall x\in F\tag 1$$ Using the fact that the inclusion $Q\subset N$ is rigid, (1) implies that $$||E_N(\theta(u))-u||_2=||\phi(u)-u||_2\leq 1/4,\forall u\in\Cal U(Q) \tag 2$$
Since $Q\subset N$, (2) implies that $$||\theta(u)u^*-1||_2^2=2-2\Re\tau(\theta(u)u^*)=\tag 3$$ $$2-2\Re\tau(E_N(\theta(u))u^*)=2\Re\tau((u-E_N(\theta(u)))u^*)\leq$$ $$ 2||u-E_N(\theta(u))||_2\leq 1/2,\forall u\in\Cal U(Q).$$
Next, we use a standard averaging trick.
For this, let $K$ denote the $||.||_2$-closure of the convex hull of the set $\{\theta(u)u^*|u\in \Cal U(Q)\}$ and let    $\xi\in K$ be the element of minimal norm. Using (3) and the fact that $K\subset (M_j)_1$, we deduce that $||\xi||\leq 1$ and that $||\xi-1||_2\leq 1/2$, so, in particular, $\xi\not=0$. Moreover, since $K$ is invariant under the $||.||_2$--preserving transformations $K\ni\eta\rightarrow\theta(u)\eta u^*$, for all $u\in\Cal U(Q)$,  the uniqueness of $\xi$ implies that $\theta(u)\xi u^*=\xi,$ for all $u\in\Cal U(Q).$ Furthermore, it is easy to see that this relation is still verified if we replace $\xi$ by the partial isometry $v$ in its polar decomposition and the unitary $u\in Q$ by an arbitrary element $x\in Q$ (since any element in C$^*$-algebra is a linear combination of 4 unitaries), i.e. $$\theta(x)v=vx,\forall x\in Q\tag 4$$ 
Using (4) it follows immediately that $v^*v\in Q'\cap M_j$ and that $vv^*\in\theta(Q)'\cap M_j=\theta(Q'\cap M_i)$. 
 Denote $q=vv^*,p_i=\theta^{-1}(q)$ and $p_j=v^*v$. Since $Q'\cap M_k=A_k$, for every $k\in I$, we get that $p_i\in A_i$ and $p_j\in A_j$.  

Now, if we define $\delta(x)=v^*\theta(x)v$, for all $x\in p_iM_ip_i$, then $\delta:p_iM_ip_i\rightarrow p_jM_jp_j$ is an isomorphism. Moreover, (4) implies that for all $x\in Q$ we have that $$\delta(xp_i)=v^*\theta(xp_i)v=v^*\theta(x)v=v^*vx=xp_j\tag 5$$

In particular, (5) implies that $\delta((Qp_i)'\cap p_iM_ip_i)=(Qp_j)'\cap p_jM_jp_j$, or, equivalently, that $\delta(A_ip_i)=A_jp_j$. Altogether, we get that $\delta$ gives an isomorphism of the inclusions $(Ap_i\subset p_iMp_i)\simeq (Ap_j\subset p_jMp_j)$. Finally, since $p_i$ and $p_j$  have the same trace, we would get that $(A_i\subset M_i)\simeq (A_j\subset M_j)$, i.e. the actions $\sigma_i$ and $\sigma_j$ are orbit equivalent ([FM]), a contradiction.\hfill$\square$
\vskip 0.1in

{\it Acknowledgments.} I am grateful to Professors Damien Gaboriau, Alekos Kechris and Sorin Popa for useful discussions and helpful comments.

\head References. \endhead
\item {[BG]} S. I. Bezuglyi, V. Ya. Golodets: {\it Hyperfinite and II1 actions for nonamenable groups}, J. Funct. Anal.  {\bf 40} (1981), no. 1, 30--44 
\item {[Bu]} M. Burger: {\it Kazhdan constants for SL$(3,\Bbb Z)$}, J. Reine Angew. Math. {\bf 413} (1991), 36--67.
\item {[C1]} A. Connes: {\it Classification of injective factors. Cases $II\sb{1},$ $II\sb{\infty },$ $III\sb{\lambda },$ $\lambda \not=1$}.  Ann. of Math. (2)  {\bf 104}(1976), no. 1, 73--115.
\item {[C2]} A. Connes: {\it A factor of type II$_1$ with countable fundamental group.}  J. Operator Theory  
{\bf 4}(1980), no. 1, 151--153.
\item {[CW]} A. Connes, B. Weiss: {\it Property (T) and asymptotically invariant sequences}, Israel J. Math. {\bf 37}(1980), 209--210.
\item {[DGRS]} H. Dooley, V. Ya. Golodets, D. J. Rudolph, D. Sinel\o{'}shchikov: {\it Non-Bernoulli systems with completely positive entropy}, Ergodic Theory Dynam. Systems, {\bf 28} (2008), 87--124.
\item {[D]} H. Dye: {\it On groups of measure preserving transformations I},  Amer. J. Math., {\bf 81} (1959), 119--159.
\item {[E]} I. Epstein: {\it Orbit inequivalent actions of non-amenable groups}, preprint math. GR 0707.4215.
\item {[Fe]} T. F\'{e}rnos: {\it Relative property (T) and linear groups}, 
Ann. Inst. Fourier (Grenoble) {\bf 56} (2006), no. 6, 1767--1804.
\item {[FM]} J. Feldman, C.C. Moore: {\it Ergodic equivalence relations, cohomology, and von Neumann algebras, II}, Trans. Am. Math. Soc. {\bf 234}(1977), 325--359. 
\item {[Fu1]} A. Furman: {\it Orbit equivalence rigidity,}  Ann. of Math. (2)  {\bf 150}(1999),  no. 3, 1083--1108.
\item {[Fu2]} A. Furman: {\it On Popa's cocycle superrigidity theorem},  IMRN  2007,  no. 19, Art. ID rnm073, 46 pp.
\item {[G]} D. Gaboriau: {\it Examples of Groups that are Measure Equivalent to the Free Group}, Ergodic Theory Dynam. Systems {\bf 25} (2005), no.6, 1809--1827. 
\item {[GG]} S.L. Gefter, V. Y. Golodets: {\it Fundamental groups for ergodic actions and actions with unit fundamental groups},  Publ. Res. Inst. Math. Sci.  {\bf 24} (1988), no. 6, 821--847.
\item {[GL]} D. Gaboriau, R. Lyons: {\it A Measurable-Group-Theoretic Solution to von Neumann's Problem}, preprint math. GR 0711.1643. 
\item {[GP]} D. Gaboriau, S.Popa: {\it An uncountable family of non orbit equivalent actions of $\Bbb F_n$}, J. Am. Math. Soc. {\bf 18}(2005), 547--559.
\item{[H]} G. Hjorth: {\it A converse to Dye's theorem}, Trans. Am. Math. Soc. {\bf 357}(2005), 3083--3103.
\item {[I1]} A. Ioana: {\it A relative version of Connes' $\chi(M)$ invariant and existence of orbit inequivalent actions},  Ergodic Theory Dynam. Systems  {\bf 27}(2007),  no. 4, 1199--1213.
\item {[I2]} A. Ioana: {\it A Construction of Non-Orbit Equivalent Actions for $\Bbb F_n$,} preprint math. GR/0610452.
\item {[K]} D. Kazhdan: {\it On the connection of the dual space of a group with the structure of its closed subgroups}, Funct. Anal. and its Appl. {\bf 1}(1967), 63--65.
\item {[Ke]} A. Kechris: {\it Global aspects of ergodic group actions and equivalence relations,} preprint 2006.
\item {[Ki]} Y. Kida: {\it Classification of certain generalized Bernoulli actions of mapping class groups,} preprint 2007.
\item {[L]} W. L\''{u}ck: {\it The type of the classifying space for a family of subgroups},  J. Pure Appl. Algebra  {\bf 149}(2000), no. 2, 177--203.
\item {[M]} G. Margulis: {\it Finitely-additive invariant measures on Euclidian spaces,}  Ergodic Theory Dynam. Systems {\bf 2}(1982), 383--396.
\item {[MSh]} N. Monod, Y. Shalom: {\it  Orbit equivalence rigidity and bounded cohomology} Ann. Math. {\bf 164} (2006) 825-–878.
\item{[MvN]} F. Murray, J. von Neumann: {\it Rings of operators}, IV, Ann. Math. {\bf 44}(1943), 716--808.
\item {[O]} A. Yu. Ol'shanskii : {\it On the question of the existence of an invariant mean on a group,} Uspekhi
Mat. Nauk, {\bf 35}(1980), no.4(214), 199–-200.
\item {[OW]} D. Ornstein, B.Weiss: {\it Ergodic theory of amenable groups. I. The Rokhlin lemma.}, Bull. Amer. Math. Soc. (N.S.) {\bf 1} (1980), 161--164.
\item {[P1]} S. Popa: {\it On a class of type II$_1$ factors with Betti numbers invariants}, Ann. of Math. {\bf 163}(2006), 809--889.
\item {[P2]} S. Popa: {\it Some computations of 1-cohomology groups and construction of non-orbit-equivalent actions}, Journal of the Inst. of Math. Jussieu {\bf 5} (2006), 309--332. 
\item {[P3]} S. Popa: {\it Strong Rigidity of II$_1$ Factors Arising from Malleable Actions of w-Rigid Groups II}, Invent. Math. {\bf 165} (2006), 409--451.
\item {[P4]} S. Popa: {\it On the superrigidity of malleable actions with spectral gap}, J. Amer. Math. Soc., to appear.  
\item {[Sh1]} Y. Shalom: {\it  Bounded generation and Kazhdan's property (T)},  Inst. Hautes Études Sci. Publ. Math.  {\bf90}(1999), 145--168. 
\item {[Sh2]} Y. Shalom: {\it  Measurable group theory,} European Congress of Mathematics,  391--423, Eur. Math. Soc., Zürich, 2005.
\item {[Sc]} K. Schmidt: {\it Amenability, Kazhdan's property T, strong ergodicity and invariant means for ergodic
group-actions}, Ergodic Theory Dynam. Systems {\bf 1}(1981), no. 2, 223–-236.
\item {[Sz]} R. Szwarc: {\it An analytic series of irreducible representations of the free group}, Annales de l'institut Fourier {\bf 38} (1988), 87--110.
\item {[T]]} A. T\H{o}rnquist: {\it Orbit equivalence and actions of $\Bbb F\sb n$},  J. Symbolic Logic  {\bf 71}(2006),  no. 1, 265--282. 
\item {[Z]} R. Zimmer: {\it Ergodic Theory and Semisimple Groups}, Birkhauser, Boston, 1984.
\enddocument